\RequirePackage{fix-cm}
\documentclass[smallextended]{svjour3}       
\smartqed  
\usepackage{graphicx}
\usepackage{fullpage}
\usepackage{subfigure}
\usepackage{amssymb}
\begin{document}

\title{Numerical solution of fractional Bratu type equations with Legendre reproducing kernel method
}\author{Mehmet Giyas Sakar\and Onur Sald{\i}r \and Ali Akg\"{u}l }


\institute{Mehmet Giyas Sakar \at
Yuzuncu Yil University, Faculty of Sciences, Department of
Mathematics, 65080, Van, Turkey \\
Tel.: +90 (432) 2251701\\
\email{giyassakar@hotmail.com}           
\and
Onur Sald{\i}r \at
Yuzuncu Yil University, Faculty of Sciences, Department of
Mathematics, 65080, Van, Turkey\\
\email{onursaldir@gmail.com}           
\and
Ali Akg\"{u}l \at
Mathematics, 56100, Siirt, Turkey\\
\email{aliakgul00727@gmail.com}           
}

\date{Received: date / Accepted: date}
\maketitle

\begin{abstract}
In this research, a new numerical method is proposed for solving fractional Bratu type boundary value problems. Fractional derivatives are taken in Caputo sense. This method is predicated on iterative approach of reproducing kernel Hilbert space theory with shifted Legendre polynomials. Construction and convergence of iterative process are shown by orthogonal projection operator. Numerical results show that our method is effective and convenient for fractional Bratu type problem.
\end{abstract}

\section{Introduction}

In recent years, a great deal of implementation of nonlinear boundary value problems with fractional order derivative take place in many areas of sciences and engineering. Various phenomena can be successfully modelled with the aid of fractional differential equations [1]. Applications, methods to find approximate solution and qualitative behaviors of solution for fractional differential equation have
been investigated by authors [2,3,4] and references therein.

In this article, a new method which is called Legendre -reproducing kernel method (L-RKM) will be proposed for attain the
approximate solution of following nonlinear fractional Bratu differential equation with Caputo derivative [5],
\begin{eqnarray}
{\ }^{c}D^{\alpha }y(x)+\lambda e^{y(x)} =0,\quad x\in\lbrack 0,1]
\end{eqnarray}
with boundary conditions
\begin{eqnarray}
y(0)=\gamma_{1},\,\ y(1)=\gamma_{2},\,\ 1<\alpha\leq2.
\end{eqnarray}
Here, $y(x)\in$ $W_{\rho}^{m}[0,1]$. Without loss of generality, we just take account of
boundary conditions $y(0)=0$ and $y(0)=0$, for $y(0)=\beta_{1}$ and
$y(1)=\beta_{2}$ can be easily reduced to $y(0)=0$ and $y(1)=0$. The Bratu's problem has no solution, one or two solutions when $\lambda>\lambda_c$, $\lambda=\lambda_c$ and $\lambda<\lambda_c$ respectively, where the critical value $\lambda_c$ is given by $\lambda_{c}=3,513830719$ for $\alpha=2$.

There are many studies about analytical and numerical solution of Bratu's problem. B-spline method [6],  Laplace decomposition method [7], Lie-group shooting method [8], homotopy analysis method [9], weighted residual method [10], decomposition method [11], perturbation iteration method [12].

The history of reproducing kernel concept goes back to the research of Zaremba in 1908
[13]. His work mainly focused on boundary value problems with Dirichlet condition
about harmonic and biharmonic functions. The reproducing property of kernel
functions has an important role in reproducing kernel Hilbert space theory. The solution of differential equations is given by convergent series form in reproducing kernel method. Recently, reproducing
kernel method is applied for numerous type of differential and integral
equations. For instance, nonlinear fourth order multi point boundary value problems [14],
Riccati differential equations [15], singular integral equations [16], fractional Bratu-type equations [17], variable order fractional functional differential equations [18], singularly perturbed problems [19], fractional advection-dispersion equation [20], fractional Riccati equation [21] and [22,23,24,25].

This paper is arranged as follows: some fundamental definitions of fractional
calculus are introduced in Section 2. Bases function and polynomial
reproducing kernel function are given in Section 3. The structure solution of
the problem with L-RKM is given in Section 4. Section 5 is reserved to
numerical examples which are successfully solved by L-RKM algorithm. Finally,
Section 6 ends with a brief conclusion.

\section{Fractional calculus}

In this section we provide some important definitions which will be used in
this study. \newline\newline\textbf{Definition 2.1} Let $x\in[0,X]$ and
$y(x)\in C[0,X]$. Then, the Riemann-Liouville fractional integral operator of
order $\alpha$ is defined as [26]:
\begin{eqnarray*}
J_{0+}^{\alpha}y(x)=\frac{1}{\Gamma(\alpha)}\int\limits_{0}^{x} {(x-r)^{\alpha
-1}y(r)dr},
\end{eqnarray*}
where $\Gamma(.)$ is Gamma function, $\alpha\geq0$ and $x>0$. Properties of
the operator $J^{\alpha}$ can be found in [1-3] and we mention only the
following: For $\alpha,\beta\ge0,$ and $\gamma> - 1:$
\begin{eqnarray*}
&&i. \,\,\ J^{0}_{0+}y(x) = y(x)\\
&&ii. \,\,\ J^{\alpha}_{0+} J^{\beta}_{0+} y(x) = J^{\alpha+ \beta}%
_{0+}y(x)\\
&&iii. \,\,\ J^{\alpha}_{0+} J^{\beta}_{0+} y(x) = J^{\beta}_{0+} J^{\alpha
}_{0+} y(x)\\
&&iv. \,\,\ J^{\alpha}_{0+} x^{\gamma}= \frac{\Gamma(\gamma+ 1)}%
{\Gamma(\alpha+ \gamma+ 1)}x^{\alpha+ \gamma}%
\end{eqnarray*}
\textbf{Definition 2.2} Let $x\in[0,X]$ and $y(x)\in AC[0,X]$. Then, the
Caputo derivative of order $\alpha$ is defined as [26]:
\begin{eqnarray*}
{\ }^{c}D_{0+}^{\alpha}y(x)=\frac{1}{\Gamma(n-\alpha)}\int_{0}^{x}%
\frac{\partial^{n}}{\partial x^{n}}\frac{y(r)}{(x-r)^{n-\alpha}}dr
\end{eqnarray*}
here $n-1<\alpha< n$, $n\in N$ and $x>0$.
\newline Also, we need
following basic properties:
\begin{eqnarray*}
&&i. \,\,\ J_{0+}^{\alpha}{\ }^{c}D_{0+}^{\alpha}y(x)=y(x)-\sum_{k=0}%
^{n-1}y^{k}(0^{+})\frac{x^{k}}{k!}\\
&&ii. \,\,\ {\ }^{c}D_{0+}^{\alpha}J_{0+}^{\alpha}y(x)=y(x)\\
&&iii. \,\,\ {\ }^{c}D_{0+}^{\alpha} y(x) = J^{n - \alpha}_{0+}D^{n} y(x)
\end{eqnarray*}

\section{Basis functions and Legendre reproducing kernel}

\subsection{\textbf{Basis functions}}

The first kind of shifted Legendre polynomials are
described on $[0,1]$ and can be obtained by the following
iterative formula:
\begin{eqnarray*}
&&P_{0}(x)  =1,\\
&&P_{1}(x)  =2x-1,\\
&&  \vdots\\
&&(n+1)P_{n+1}(x)=(2n+1)(2x-1)P_{n}(x)-nP_{n-1}(x),
\end{eqnarray*}
for $n=1,2,...$. The orthogonality requirement is
\begin{eqnarray}
\langle P_{n},P_{m} \rangle&=&\int_{0}^{1}\rho_{[0,1]}(x)P_{n}(x)P_{m}(x)dx\nonumber\\
&=&\left\{
\begin{array}
[c]{ll}%
0, & n\neq m,\\
1, & n=m=0,\\
\frac{1}{2n+1}, & n=m\neq0,
\end{array}
\right.
\end{eqnarray}
where weighted function is taken as,
\begin{eqnarray}
\rho_{[0,1]}(x)=1.
\end{eqnarray}
We can establish Legendre basis
functions which provide the homogeneous boundary conditions as:
\begin{eqnarray}
y(0)=0 \,\,\hbox{and}\,\ y(1)=0.
\end{eqnarray}
This is useful property for solving boundary value problems. Now, the basis functions defined as
\begin{eqnarray}
\varphi_{i}(x)= \left\{
\begin{array}
[c]{ll}%
P_{i}(x)-P_{0}(x), & \hbox{ $i$ is even,}\\
P_{i}(x)-P_{1}(x), & \hbox{$i$ is odd.}
\end{array}
\right.  i\geq2.
\end{eqnarray}
have the conditions
\begin{eqnarray}
\varphi_{i}(0)=\varphi_{i}(1)=0.
\end{eqnarray}
It is well known that from the theory of Legendre polynomials, this basis functions constructed by (6) are complete system.

\subsection{\textbf{Legendre reproducing kernel function}}

\noindent\textbf{Definition 3.1} Let $X$ be a nonempty set, and $H$ be a
Hilbert space of real-valued functions on some set $X$ with inner product
$\langle\cdot,\cdot\rangle_{H}$. Then, function $K:X\times X\rightarrow R$ is
called to be the reproducing kernel function of $H$ if and only if

\begin{enumerate}
\item $K(x,\cdot) \in H, \forall x \in X$

\item $\langle\varphi(\cdot),K(x,\cdot) \rangle= \varphi(x), \forall\varphi\in
H, \forall x\in X$.
\end{enumerate}

The second condition is reproducing kernel property. Then, this Hilbert space is called reproducing kernel Hilbert space (RKHS).
\newline\newline\textbf{Theorem 3.1} Let $H$ be $n$-dimensional Hilbert space,
$\{e_{i}\}_{i=1}^{n}$ is an orthonormal basis of $H$, then the reproducing kernel of $H$ as [27]:
\begin{eqnarray}
K_{t}(x)=\sum_{i=1}^{n}e_{i}(x)\bar{e}_{i}(t).
\end{eqnarray}
\textbf{Definition 3.2}
Let $W_{\rho}^{m}[0,1]$ be the weighted inner
product space of polynomials on $[0,1]$ with real coefficients and degree less
than or equal to $m$ with inner product
\[
\langle y,v\rangle_{W_{\rho}^{m}}=\int_{0}^{1}\rho_{\lbrack 0,1]}%
(x)y(x)v(x)dx,\,\,\ \forall y,v\in W_{\rho}^{m}[0,1]
\]
with $\rho_{\lbrack 0,1]}(x)$ defined by equation (4), and the norm
\[
\Vert y\Vert_{W_{\rho}^{m}}=\sqrt{\langle y,y\rangle}_{W_{\rho}^{m}%
},\,\,\ \forall y\in W_{\rho}^{m}[0,1].
\]
From definiton of $L_{\rho}^{2}[0,1]=\{f|\int_{0}^{1}\rho_{\lbrack 0,1]}(x)|f(x)|^{2}%
dx<\infty\}$ for any fixed $m$, $W_{\rho}%
^{m}[0,1]$ is a subspace of $L_{\rho}^{2}[0,1]$ and $\forall y,v\in W_{\rho
}^{m}[0,1]$, $\langle y,v\rangle_{W_{\rho}^{m}}=\langle y,v\rangle_{L_{\rho
}^{2}}$
\newline\newline\textbf{Theorem 3.2} $W_{\rho}%
^{m}[0,1]$ function space with its inner product and norm is a RKHS.
\newline\newline\textbf{Proof.} From Definition 3.2, $W_{\rho}^{m}[0,1]$ is
a finite dimensional inner product space. Every finite-dimensional inner product space is a Hilbert space. Therefore, from
this result and Theorem 3.1, $W_{\rho}^{m}[0,1]$ is a RKHS.\newline\newline For solving boundary value problem,
it is necessary to define a closed subspace of $W_{\rho}^{m}[0,1]$ which satisfy homogenous boundary conditions.
\newline\newline\textbf{Definition 3.3} Let
\begin{eqnarray*}
{\ }^{0}W_{\rho}^{m}[0,1]=\{y|y\in W_{\rho}^{m}[0,1],y(0)=y(1)=0\}.
\end{eqnarray*}
By using equation (6), we can easily show
that ${\ }^{0}W_{\rho}^{m}[0,1]$ function space is a RKHS. From Theorem 3.1, the kernel
function $K_{t}^{m}(x)$ of ${\ }^{0}W_{\rho}^{m}[0,1]$ can be written as follow
\begin{eqnarray}
K_{t}^{m}(x)=\sum_{i=2}^{m}h_{i}(x){h_{i}}(t).
\end{eqnarray}
Here, $h_{i}(x)$ is complete orthonormal system which is derived from Eq. (6)
by use of Gram-Schmidt process. Eq. (9) is very practical for application. Namely, $K_{t}^{m}(x)$ and $W_{\rho}^{m}[0,1]$ space can easily re-computed by ascending $m$.

\section{The Legendre reproducing kernel method (L-RKM)}

\subsection{Representation of exact solution in ${\ }^{0}W_{\rho}^{m}[0,1]$}

In this subsection, we construct a Legendre reproducing kernel method for analytical solution of Eqs. (1)-(2). The solution of (1)-(2) will be given in ${\ }%
^{0}W_{\rho}^{m}[0,1]$. First of all, the linear operator $L$ is defined as,,
\[
L:{\ }^{0}W_{\rho}^{m}[0,1]\rightarrow L_{\rho}^{2}[0,1]
\]
such that
\[
Ly(x):={\ }^{c}D^{\alpha}y(x)
\]
The Eq.(1)-(2) can be written as follows
\begin{eqnarray}
\left\{
\begin{array}
[c]{ll}%
Ly=-\lambda e^{y(x)} & \\
y(0)=y(1)=0. &
\end{array}
\right.
\end{eqnarray}
Let $g(x,y)$=$-\lambda \exp (y(x))$. It is easy to show that $L$ is a bounded linear operator. We shall give the
representation of an analytical solution of equation (10) in the space
${\ }^{0}W_{\rho}^{m}[0,1]$. Let $K_{t}^{m}(x)$ be the polynomial reproducing
kernel function of ${\ }^{0}W_{\rho}^{m}[0,1]$.\newline\newline\textbf{Theorem
4.1} For Eqs. (1)-(2), if $\{x_{i}\}_{i=0}^{m-2}$ be any $(m-1)$ distinct
points in $(0,1)$, then $\psi_{i}^{m}(x)=L^{\ast}K_{x_{i}}^{m}(x)=L_{t}%
K_{t}^{m}(x)|_{t=x_{i}}$.\newline\newline\textbf{Proof} For any fixed
$x_{i}\in(0,1)$, put
\begin{eqnarray}
\psi_{i}^{m}(x)  &  =L^{\ast}K_{x_{i}}^{m}(x)=\langle L^{\ast}K_{x_{i}}%
^{m}(x),K_{x}^{m}(t)\rangle_{{\ }^{0}W_{\rho}^{m}}\nonumber\\
&  =\langle K_{x_{i}}^{m}(x),L_{t}K_{x}^{m}(t)\rangle_{L_{\rho}^{2}}%
=L_{t}K_{x}^{m}(t)|_{t=x_{i}}.
\end{eqnarray}
From Definition 3.3, it is clear that $K_{x}^{m}(t)=K_{t}^{m}(x)$. Therefore
$\psi_{i}^{m}(x)=L^{\ast}K_{x_{i}}^{m}(x)=L_{t}K_{t}^{m}(x)|_{t=x_{i}}$. Here
$L^{\ast}$ is the adjoint operator of $L$ and the subscript $t$ of $L$ shows that the operator $L$ applies to the function $t$. For any fixed $m$ and $x_{i}\in(0,1)$, $\psi_{i}^{m}\in{\ }%
^{0}W_{\rho}^{m}[0,1]$.\newline\newline\textbf{Theorem 4.2} For $m\geq2$, let
$\{x_{i}\}_{i=0}^{m-2}$ be any $(m-1)$-distinct points in $(0,1)$, then
$\{\psi_{i}^{m}\}_{i=0}^{m-2}$ is a complete system of ${\ }^{0}W_{\rho}%
^{m}[0,1]$.\newline\newline\textbf{Proof.} For each fixed $u\in{\ }^{0}%
W_{\rho}^{m}[0,1]$, let
\[
\langle y(x),\psi_{i}^{m}(x)\rangle_{{\ }^{0}W_{\rho}^{m}}%
=0,\,\ i=0,1,...,m-2,
\]
this means, for $i=0,1,...,m-2$,
\begin{eqnarray}
\langle y(x),\psi_{i}^{m}(x)\rangle_{{\ }^{0}W_{\rho}^{m}}  &=&\langle
y(x),L^{\ast}K_{x_{i}}^{m}(x)\rangle_{{\ }^{0}W_{\rho}^{m}}\nonumber\\
&=&\langle Ly(x),K_{x_{i}}^{m}(x)\rangle_{L_{\rho}^{2}}\nonumber\\
&=&Ly(x_{i})=0
\end{eqnarray}
By the existence of inverse operator $L^{-1}$ it is concluded that
$y\equiv0$. Therefore, $\{\psi_{i}^{m}\}_{i=0}^{m-2}$ is a complete system of
${\ }^{0}W_{\rho}^{m}[0,1]$. Proof is completed.\newline\newline Theorem 4.2
shows that in L-RKM, use of a finite distinct points is sufficient and does not need a
dense sequence. So, our approach is differ from traditional reproducing kernel method.\\
The orthonormal system $\{\bar{\psi}_{i}%
^{m}\}_{i=0}^{m-2}$ of ${\ }^{0}W_{\rho}^{m}[0,1]$ can be obtained from the
Gram-Schmidt orthogonalization process using $\{\psi_{i}^{m}%
\}_{i=0}^{m-2}$,
\begin{eqnarray}
\bar{\psi}_{i}^{m}(x)=\sum_{k=0}^{i}\beta_{ik}^{m}\psi_{k}^{m}(x),
\end{eqnarray}
where $\beta_{ik}^{m}$ are orthogonalization coefficients.\newline%
\newline\textbf{Theorem 4.3} Assume that $y_{m}$ is the exact solution
of Eqs. (1)-(2) in ${\ }^{0}W_{\rho}^{m}[0,1]$. Let $\{x_{i}\}_{i=0}%
^{m-2}$ be any $(m-1)$ distinct points in $(0,1)$, then
\begin{eqnarray}
y_{m}(x)=\sum_{i=0}^{m-2}\sum_{k=0}^{i}\beta_{ik}^{m}g(x_{k},y_{m}(x_{k}%
))\bar{\psi}_{i}^{m}(x).
\end{eqnarray}
\textbf{Proof.} Since $y_{m}\in{\ }^{0}W_{\rho}^{m}[0,1]$ by Theorem 4.2 we
can write
\[
y_{m}(x)=\sum_{i=0}^{m-2}\langle y_{m}(x),\bar{\psi}_{i}^{m}(x)\rangle
_{{\ }^{0}W_{\rho}^{m}}\bar{\psi}_{i}^{m}(x).
\]
Otherwise, by using Eq. (11) and Eq. (13), we obtain $y_{m}(x)$ which the exact solution of Eq. (10) in ${\ }^{0}W_{\rho}^{m}[0,1]$ as,
\begin{eqnarray*}
y_{m}(x) &=&\sum_{i=0}^{m-2}\langle y_{m}(x),\bar{\psi}_{i}^{m}%
(x)\rangle_{{\ }^{0}W_{\rho}^{m}}\bar{\psi}_{i}^{m}(x)\\
&=&\sum_{i=0}^{m-2}\langle y_{m}(x),\sum_{k=0}^{i}\beta_{ik}^{m}\psi_{k}%
^{m}(x)\rangle_{{\ }^{0}W_{\rho}^{m}}\bar{\psi}_{i}^{m}(x)\\
&=&\sum_{i=0}^{m-2}\sum_{k=0}^{i}\beta_{ik}^{m}\langle y_{m}(x),\psi_{k}%
^{m}(x)\rangle_{{\ }^{0}W_{\rho}^{m}}\bar{\psi}_{i}^{m}(x)\\
&=&\sum_{i=0}^{m-2}\sum_{k=0}^{i}\beta_{ik}^{m}\langle y_{m}(x),L^{\ast
}K_{x_{k}}^{m}(x)\rangle_{{\ }^{0}W_{\rho}^{m}}\bar{\psi}_{i}^{m}(x)\\
&=&\sum_{i=0}^{m-2}\sum_{k=0}^{i}\beta_{ik}^{m}\langle Ly_{m}(x),K_{x_{k}%
}^{m}(x)\rangle_{L_{\rho}^{2}}\bar{\psi}_{i}^{m}(x)\\
&=&\sum_{i=0}^{m-2}\sum_{k=0}^{i}\beta_{ik}^{m}\langle g(x,y_{m}%
(x)),K_{x_{k}}^{m}(x)\rangle_{L_{\rho}^{2}}\bar{\psi}_{i}^{m}(x)\\
&=&\sum_{i=0}^{m-2}\sum_{k=0}^{i}\beta_{ik}^{m}g(x_{k},y_{m}(x_{k}))\bar
{\psi}_{i}^{m}(x).
\end{eqnarray*}
The proof is completed. \newline\newline\textbf{Theorem 4.4} If
$y_{m}(x)\in{\ }^{0}W_{\rho}^{m}[0,1]$, then $|y_{m}^{(r)}(x)|\leq F\Vert
y_{m}\Vert_{{\ }^{0}W_{\rho}^{m}}$ for $r=0,\ldots,m-1$, where $F$ is a
constant. \newline\newline\textbf{Proof.} For any $x,\,t\in\left[
{0,1}\right]  $, we have
\begin{eqnarray*}
y_{m}^{(r)}\left(  x\right)  =\langle y_{m}\left(
t\right)  ,\partial_{x}^{r}K_{x}\left(  t\right)  \rangle_{{\ }^{0}W_{\rho
}^{m}},
\end{eqnarray*}
$r=0,\ldots,m-1.$ By the expression of $K_{x}\left(  t\right)  $, it
follow that $\left\Vert {\partial_{x}^{r}K_{x}\left(  t\right)  }\right\Vert
_{{\ }^{0}W_{\rho}^{m}}\leq F_{r},\,r=0,\ldots,m-1.$\newline So,
\begin{eqnarray*}
|y_{m}^{(r)}(x)|  & =&|{\langle y_{m}(x),\partial_{x}^{r}K_{x}(t)\rangle
_{{\ }^{0}W_{\rho}^{m}}}|\\
&\leq&\Vert{y_{m}(x)}\Vert_{{\ }^{0}W_{\rho}^{m}[a,b]}\Vert{\partial_{x}%
^{r}K_{x}(x)}\Vert_{{\ }^{0}W_{\rho}^{m}}\\
&\leq& F_{r}\Vert{y_{m}(x)}\Vert_{{\ }^{0}W_{\rho}^{m}},r=0,\ldots,m-1.
\end{eqnarray*}
Therefore, $|y_{m}^{(r)}(x)|\leq\max\{F_{0},\ldots,F_{m-1}\}\left\Vert
{y_{m}\left(  x\right)  }\right\Vert _{{\ }^{0}W_{\rho}^{m}}$, $r=0,\ldots
,m-1$. This completes the proof. \newline\newline\textbf{Theorem 4.5} The
approximate solution $y_{m}(x)$ and its derivatives $y_{m}^{(r)}(x)$ are
uniformly converge to $y(x)$ and $y^{(r)}(x)$ ($r=0,\ldots,m-1$) respectively.
\newline\newline\textbf{Proof} From Theorem 4.4 for any $x\in\lbrack 0,1]$ we
get
\begin{eqnarray*}
|y_{m}^{(r)}(x)-y^{(r)}(x)| &=&|\langle y_{m}(x)-y(x),\partial_{x}^{r}%
K_{x}(x)\rangle|_{{\ }^{0}W_{\rho}^{m}}\\
&\leq&\Vert\partial_{x}^{r}K_{x}(x)\Vert_{{\ }^{0}W_{\rho}^{m}}\Vert
y_{m}(x)-y(x)\Vert_{{\ }^{0}W_{\rho}^{m}}\\
&\leq& F_{r}\Vert y_{m}(x)-y(x)\Vert_{{\ }^{0}W_{\rho}^{m}},\,\ r=0,\ldots
,m-1.
\end{eqnarray*}
where $F_{0},\ldots,F_{m-1}$ are positive constants. Therefore, if $y_{m}%
(x)\rightarrow y(x)$ in the sense of the norm of ${\ }^{0}W_{\rho}^{m}[0,1]$
as $m\rightarrow\infty$, the approximate solution $y_{m}(x)$ and its
derivatives $y_{m}^{^{\prime}}(x),...,y_{m}^{(m-1)}(x)$ are uniformly
converge to the exact solution $y(x)$ and its derivatives $y^{^{\prime}%
}(x),...,y^{(m-1)}(x)$ respectively.\newline\newline\textbf{Remark 4.1} We
will pay attention the next two states in order to solve Eqs. (1)-(2) by employing
L-RKM. \newline\textbf{Case. 1}: The approximate
solution can be attained directly from Eq. (14), if Eq. (1) is linear. \newline\textbf{Case. 2}: The approximate solution can be attained by using the
following iterative process, if Eq.(1) is nonlinear.

\subsection{\textbf{Construction and convergence of iterative process}}

\noindent Firstly, we constitute the following iterative sequence $y_{m}(x)$,
putting,
\begin{equation}
\left\{
{{\begin{array}{*{20}c} {Lu_{m,n}\left( x \right) = g\left( {x,y_{m,n-1}(x)} \right)} \hfill \\ {y_{m,n}\left( x \right) = P_{m-1} u_{m,n} (x)} \hfill \\ \end{array}}%
}\right.
\end{equation}
where $P_{m-1}:{\ }^{0}W_{\rho}^{m}[0,1]\rightarrow span\{\bar{\psi}_{0}%
^{m},\bar{\psi}_{1}^{m},\ldots,\bar{\psi}_{m-2}^{m}\}$ is a orthogonal
projection operator and $u_{m,n}(x)\in{\ }^{0}W_{\rho}^{m}[0,1]$ is the $n$-th
iterative approximate solution of (15). Secondary, we will give a Lemma for
construction of iterative process. \newline\newline\textbf{Theorem 4.6} If
$\{x_{i}\}_{i=0}^{m-2}$ is distinct points in $(0,1)$, then
\[
u_{m,n}(x)=\sum_{i=0}^{m-2}\sum_{k=0}^{i}\beta_{ik}^{m}g(x_{k},y_{m,n-1}%
(x_{k}))\bar{\psi}_{i}^{m}(x)
\]
\textbf{Proof.} Since $u_{m,n}(x)\in{\ }^{0}W_{\rho}^{m}[0,1]$, $\{\bar{\psi
}_{i}^{m}(x)\}_{i=0}^{m-2}$ is the complete system in ${\ }^{0}W_{\rho}%
^{m}[0,1]$,
\begin{eqnarray*}
u_{m,n}(x) &=&\sum_{i=0}^{m-2}\langle u_{m,n}(x),\bar{\psi}_{i}%
^{m}(x)\rangle_{{\ }^{0}W_{\rho}^{m}}\bar{\psi}_{i}^{m}(x)\\
&=&\sum_{i=0}^{m-2}\langle u_{m,n}(x),\sum_{k=0}^{i}\beta_{ik}^{m}\psi
_{k}^{m}(x)\rangle_{{\ }^{0}W_{\rho}^{m}}\bar{\psi}_{i}^{m}(x)\\
&=&\sum_{i=0}^{m-2}\sum_{k=0}^{i}\beta_{ik}^{m}\langle u_{m,n}(x),\psi
_{k}^{m}(x)\rangle_{{\ }^{0}W_{\rho}^{m}}\bar{\psi}_{i}^{m}(x)\\
&=&\sum_{i=0}^{m-2}\sum_{k=0}^{i}\beta_{ik}^{m}\langle u_{m,n}(x),L^{\ast
}K_{x_{k}}^{m}(x)\rangle_{{\ }^{0}W_{\rho}^{m}}\bar{\psi}_{i}^{m}(x)\\
&=&\sum_{i=0}^{m-2}\sum_{k=0}^{i}\beta_{ik}^{m}\langle Lu_{m,n}(x),K_{x_{k}%
}^{m}(x)\rangle_{L_{\rho}^{2}}\bar{\psi}_{i}^{m}(x)\\
&=&\sum_{i=0}^{m-2}\sum_{k=0}^{i}\beta_{ik}^{m}\langle g(x,y_{m,n-1}%
(x)),K_{x_{k}}^{m}(x)\rangle_{L_{\rho}^{2}}\bar{\psi}_{i}^{m}(x)\\
&=&\sum_{i=0}^{m-2}\sum_{k=0}^{i}\beta_{ik}^{m}g(x_{k},y_{m,n-1}(x_{k}%
))\bar{\psi}_{i}^{m}(x)
\end{eqnarray*}
This completes the proof.\newline Taking $y_{m,0}(x)=0$ (we can choose any fixed
$y_{m,0}(x)\in{\ }^{0}W_{\rho}^{m}[0,1]$) define the iterative sequence
\begin{eqnarray}
y_{m,n}(x)&=&P_{m-1}u_{m,n}(x)\nonumber\\
&=&\sum_{i=0}^{m-2}\sum_{k=0}^{i}\beta_{ik}%
^{m}g(x_{k},y_{m,n-1}(x_{k}))\bar{\psi}_{i}^{m}(x),\,\ n=1,2,...
\end{eqnarray}

\section{Numerical examples}

In this part, a numerical example is considered to show the efficiency and accuracy of the present method. Numerical result shows that L-RKM is powerful and convenient.\\\\
\noindent\textbf{Example 5.1} We consider the following  nonlinear fractional
Bratu equation with Caputo derivative
\begin{equation}
{\ }^{c}D^{\alpha }y(x)+\lambda e^{y(x)}=0,\,\,\,1<\alpha \leq 2,\,\,\ x\in
\lbrack 0,1]  \label{1.13}
\end{equation}%
\begin{equation}
y(0)=0,\,\,\ y(1)=0.  \label{1.14}
\end{equation}%
The exact solution for $\alpha =2$ is given by
\begin{eqnarray}
y(x)=-2\ln \left[\frac{\cosh\left(\left(x-\frac{1}{2}\right) \frac{\theta }{2}\right) }{\cosh \left( \frac{\theta }{4}\right) }\right]
\end{eqnarray}
where $\theta $ is the solution of $\theta =\sqrt{2\lambda }\cosh \left(
\theta /4\right) $.
Using L-RKM for Eqs. (17)-(18), taking $x_{i}=\frac{i+0.3}{m}$, \,\ $i=0,1,\,2,...,m-2$, the numerical solution $y_{m,n}\left(  x\right)  $ is
calculated by Eq. (15). Comparison of L-RKM at some selected grid points
with exact solution are given in Table 1, Table 4 and Table 7. Absolute error of different $m$
values for this example are given in Table 2, Table 5 and Table 8. Comparison of absolute error for different methods are given in Table 3 and Table 6.  Absolute error of L-RKM solution for $y(x)$ and its derivatives are given in Fig. 1-5.

\section{Conclusion}

In this research, Legendre reproducing kernel method (L-RKM) is proposed and
successfully applied to find the numerical solution of the fractional Bratu-type boundary value problems.
Numerical results show that the present method is a powerful and effective for
solving fractional Bratu-type boundary value problems. Furthermore, convergence of present iterative method is discussed.

\newpage

\begin{table*}[ht]
	\caption{Numerical results for Example 5.1 with fractional order ($\lambda$=1, $m=20$ and $n=30$)}
	\begin{tabular}{lllll}
		\hline\noalign{\smallskip}
		$x$ & Exact Sol.       & $\alpha$=2        & $\alpha$=1.9  & $\alpha$=1.8 \\
		\noalign{\smallskip}\hline\noalign{\smallskip}
		0.1 & 0.04984679124541 & 0.04984679124541  & 0.053348029010 & 0.05640455239 \\
		0.2 & 0.08918993462882 & 0.08918993462882  & 0.094053979156 & 0.09778884381 \\
		0.3 & 0.11760909576794 & 0.11760909576794  & 0.122561203471 & 0.12577795186 \\
		0.4 & 0.13479025388418 & 0.13479025388418  & 0.139016612388 & 0.14107528890 \\
		0.5 & 0.14053921440047 & 0.14053921440047  & 0.143580804066 & 0.14424570594 \\
		0.6 & 0.13479025388418 & 0.13479025388418  & 0.136494883917 & 0.13586263427 \\
		0.7 & 0.11760909576794 & 0.11760909576794  & 0.118101177356 & 0.11655270867 \\
		0.8 & 0.08918993462882 & 0.08918993462882  & 0.088844712904 & 0.08700454024 \\
		0.9 & 0.04984679124541 & 0.04984679124541  & 0.049263971500 & 0.04796098207 \\
		\noalign{\smallskip}\hline
	\end{tabular}
\end{table*}

\begin{table*}[ht]
	\caption{Absolute error of Example 5.1 for different $m$ values with $\lambda$=1, $n=30$ and $\alpha=2$}
	\begin{tabular}{llllll}
		\hline\noalign{\smallskip}
		$x$ &$m=10$ & $m=12$ & $m=14$ & $m=16$ & $m=18$ \\
		\noalign{\smallskip}\hline\noalign{\smallskip}
		0.1 & 4.59E-10 & 1.25E-11 & 3.52E-13 & 1.01E-14 & 2.96E-16 \\
		0.2 & 8.17E-10 & 2.30E-11 & 6.62E-13 & 1.93E-14 & 5.68E-16 \\
		0.3 & 1.18E-9  & 3.34E-11 & 9.65E-13 & 2.82E-14 & 8.35E-16 \\
		0.4 & 1.53E-9  & 4.35E-11 & 1.25E-12 & 3.69E-14 & 1.09E-15 \\
		0.5 & 1.86E-9  & 5.30E-11 & 1.53E-12 & 4.51E-14 & 1.33E-15 \\
		0.6 & 2.17E-9  & 6.19E-11 & 1.79E-12 & 5.28E-14 & 1.56E-15 \\
		0.7 & 2.46E-9  & 7.02E-11 & 2.03E-12 & 5.99E-14 & 1.77E-15 \\
		0.8 & 2.72E-9  & 7.75E-11 & 2.25E-12 & 6.63E-14 & 1.96E-15 \\
		0.9 & 2.90E-9  & 8.44E-11 & 2.45E-12 & 7.22E-14 & 2.13E-15 \\
		\noalign{\smallskip}\hline
	\end{tabular}
\end{table*}

\begin{table*}[ht]
	\caption{Absolute error of L-RKM and different methods for Example 5.1 case $\lambda=1$}
	\begin{tabular}{llllll}
		\hline\noalign{\smallskip}
		$x$ &B-Spline [6] & Laplace [7] & LGSM [8] & DM [11] & L-RKM\\
		\noalign{\smallskip}\hline\noalign{\smallskip}
		0.1 & 2.97E-6  & 1.97E-6 & 7.50E-7 & 2.68E-3 & 5.53E-17 \\
		0.2 & 5.46E-6  & 3.93E-6 & 1.01E-6 & 2.02E-3 & 1.09E-16 \\
		0.3 & 7.33E-6  & 5.85E-6 & 9.04E-7 & 1.52E-4 & 1.63E-16 \\
		0.4 & 8.49E-6  & 7.70E-6 & 5.23E-7 & 2.20E-3 & 2.15E-15 \\
		0.5 & 8.89E-6  & 9.46E-6 & 5.06E-9 & 3.01E-3 & 2.64E-15 \\
		0.6 & 8.49E-6  & 1.11E-5 & 5.13E-7 & 2.20E-3 & 3.11E-15 \\
		0.7 & 7.33E-6  & 1.25E-5 & 8.94E-7 & 1.52E-4 & 3.54E-15 \\
		0.8 & 5.46E-6  & 1.34E-5 & 1.00E-6 & 2.02E-3 & 3.92E-15 \\
		0.9 & 2.97E-6  & 1.19E-5 & 7.41E-7 & 2.68E-3 & 4.27E-15 \\
		\noalign{\smallskip}\hline
	\end{tabular}
\end{table*}

\begin{table*}[ht]
	\caption{Numerical results for Example 5.1 with fractional order ($\lambda$=2, $m=20$ and $n=30$)}
	\begin{tabular}{lllll}
		\hline\noalign{\smallskip}
		$x$ & Exact Sol.       & $\alpha$=2        & $\alpha$=1.9  & $\alpha$=1.8\\
		\noalign{\smallskip}\hline\noalign{\smallskip}
		0.1 & 0.11441074326774 & 0.11441074326804  & 0.12395416270 & 0.13161211082 \\
		0.2 & 0.20641911648760 & 0.20641911648819  & 0.22060766752 & 0.23047094173 \\
		0.3 & 0.27387931182555 & 0.27387931182641  & 0.28927167294 & 0.29804837485 \\
		0.4 & 0.31508936422566 & 0.31508936422678  & 0.32910866272 & 0.33469946332 \\
		0.5 & 0.32895242134111 & 0.32895242134245  & 0.33985942541 & 0.34125713685 \\
		0.6 & 0.31508936422566 & 0.31508936422720  & 0.32203437200 & 0.31930371842 \\
		0.7 & 0.27387931182555 & 0.27387931182723  & 0.27693051224 & 0.27114363748 \\
		0.8 & 0.20641911648760 & 0.20641911648939  & 0.20652983899 & 0.19963116013 \\
		0.9 & 0.11441074326774 & 0.11441074326959  & 0.11332086854 & 0.10793583205 \\
		\noalign{\smallskip}\hline
	\end{tabular}
\end{table*}

\begin{table*}[ht]
	\caption{Absolute error of Example 5.1 for different $m$ values with $\lambda$=2, $n=30$ and $\alpha=2$ }
	\begin{tabular}{llllll}
		\hline\noalign{\smallskip}
		$x$ &$m=10$ & $m=12$ & $m=14$ & $m=16$ & $m=18$\\
		\noalign{\smallskip}\hline\noalign{\smallskip}
		0.1 & 9.65E-8 & 6.25E-9 & 4.17E-10 & 2.84E-11 & 1.96E-12 \\
		0.2 & 1.76E-7 & 1.16E-8 & 7.92E-10 & 5.44E-11 & 3.78E-12 \\
		0.3 & 2.53E-7 & 1.69E-8 & 1.14E-9  & 7.91E-11 & 5.50E-12 \\
		0.4 & 3.24E-7 & 2.16E-8 & 1.47E-9  & 1.01E-10 & 7.08E-12 \\
		0.5 & 3.86E-7 & 2.58E-8 & 1.76E-9  & 1.21E-10 & 8.47E-12 \\
		0.6 & 4.37E-7 & 2.92E-8 & 1.99E-9  & 1.38E-10 & 9.62E-12 \\
		0.7 & 4.76E-7 & 3.19E-8 & 2.18E-9  & 1.50E-10 & 1.05E-11 \\
		0.8 & 5.02E-7 & 3.37E-8 & 2.30E-9  & 1.59E-10 & 1.11E-11 \\
		0.9 & 5.09E-7 & 3.48E-8 & 2.38E-9  & 1.64E-10 & 1.14E-11 \\
		\noalign{\smallskip}\hline
	\end{tabular}
\end{table*}

\begin{table*}[ht]
	\caption{Absolute error of L-RKM and different methods for Example 5.1 case $\lambda=2$ }
	\begin{tabular}{llllll}
		\hline\noalign{\smallskip}
		$x$ &B-Spline [6] & Laplace [7] & LGSM [8] & DM [11] & L-RKM\\
		\noalign{\smallskip}\hline\noalign{\smallskip}
		0.1 & 1.71E-5  & 2.12E-3 & 4.03E-6 & 1.52E-2 & 4.08E-13 \\
		0.2 & 3.25E-5  & 4.20E-3 & 5.70E-6 & 1.46E-2 & 8.04E-13 \\
		0.3 & 4.48E-5  & 6.18E-3 & 5.22E-6 & 5.88E-3 & 1.18E-12 \\
		0.4 & 5.28E-5  & 8.00E-3 & 3.07E-6 & 3.24E-3 & 1.52E-12 \\
		0.5 & 5.26E-5  & 9.59E-3 & 1.45E-8 & 6.98E-3 & 1.83E-12 \\
		0.6 & 5.28E-5  & 1.09E-2 & 3.04E-6 & 3.24E-3 & 2.08E-12 \\
		0.7 & 4.48E-5  & 1.19E-2 & 5.19E-6 & 5.88E-3 & 2.28E-12 \\
		0.8 & 3.25E-5  & 1.23E-2 & 5.67E-6 & 1.46E-2 & 2.42E-12 \\
		0.9 & 1.71E-5  & 1.08E-2 & 4.01E-6 & 1.52E-2 & 2.50E-12 \\
		\noalign{\smallskip}\hline
	\end{tabular}
\end{table*}

\begin{table*}[ht]
	\caption{Numerical results for Example 5.1 with fractional order ($\lambda$=3, $m=20$ and $n=30$)}
	\begin{tabular}{lllll}
		\hline\noalign{\smallskip}
		$x$ & Exact Sol.       & $\alpha$=2        & $\alpha$=1.9  & $\alpha$=1.8\\
		\noalign{\smallskip}\hline\noalign{\smallskip}
		0.1 & 0.21577498476753 & 0.21577498329811  & 0.24085457687 & 0.26065428735 \\
		0.2 & 0.39431958747803 & 0.39431958460834  & 0.43526223628 & 0.46414428666 \\
		0.3 & 0.52843643955821 & 0.52843643546802  & 0.57680617280 & 0.60602026500 \\
		0.4 & 0.61182920818459 & 0.61182920316127  & 0.65968880520 & 0.68206752260 \\
		0.5 & 0.64014669604146 & 0.64014669046258  & 0.68090952720 & 0.69191607740 \\
		0.6 & 0.61182920818459 & 0.61182920247466  & 0.64127386680 & 0.63986648210 \\
		0.7 & 0.52843643955821 & 0.52843643413164  & 0.54537595346 & 0.53412579803 \\
		0.8 & 0.39431958747803 & 0.39431958268710  & 0.40065374951 & 0.38508131538 \\
		0.9 & 0.21577498476753 & 0.21577498087398  & 0.21597352866 & 0.20346054571 \\
		\noalign{\smallskip}\hline
	\end{tabular}
\end{table*}

\begin{table*}[ht]
	\caption{Absolute error of Example 5.1 for different $m$ values with $\lambda$=3, $n=30$ and $\alpha=2$}
	\begin{tabular}{llllll}
		\hline\noalign{\smallskip}
		$x$ &$m=10$ & $m=12$ & $m=14$ & $m=16$ & $m=18$\\
		\noalign{\smallskip}\hline\noalign{\smallskip}
		0.1 & 8.91E-6 & 1.14E-6 & 1.52E-7  & 1.97E-8 & 3.42E-9 \\
		0.2 & 1.67E-5 & 2.19E-6 & 2.93E-7  & 3.81E-8 & 6.63E-9 \\
		0.3 & 2.39E-5 & 3.13E-6 & 4.20E-7  & 5.49E-8 & 9.49E-9 \\
		0.4 & 2.99E-5 & 3.92E-6 & 5.27E-7  & 6.89E-8 & 1.17E-8 \\
		0.5 & 3.43E-5 & 4.50E-6 & 6.04E-7  & 7.93E-8 & 1.33E-8 \\
		0.6 & 3.67E-5 & 4.82E-6 & 6.47E-7  & 8.52E-8 & 1.40E-8 \\
		0.7 & 3.71E-5 & 4.87E-6 & 6.54E-7  & 8.67E-8 & 1.38E-8 \\
		0.8 & 3.56E-5 & 4.68E-6 & 6.28E-7  & 8.38E-8 & 1.28E-8 \\
		0.9 & 3.21E-5 & 4.29E-6 & 5.77E-7  & 7.74E-8 & 1.32E-8 \\
		\noalign{\smallskip}\hline
	\end{tabular}
\end{table*}

\begin{figure*}[htbp]
	\centerline{\includegraphics[width=2.400in,height=1.55in]{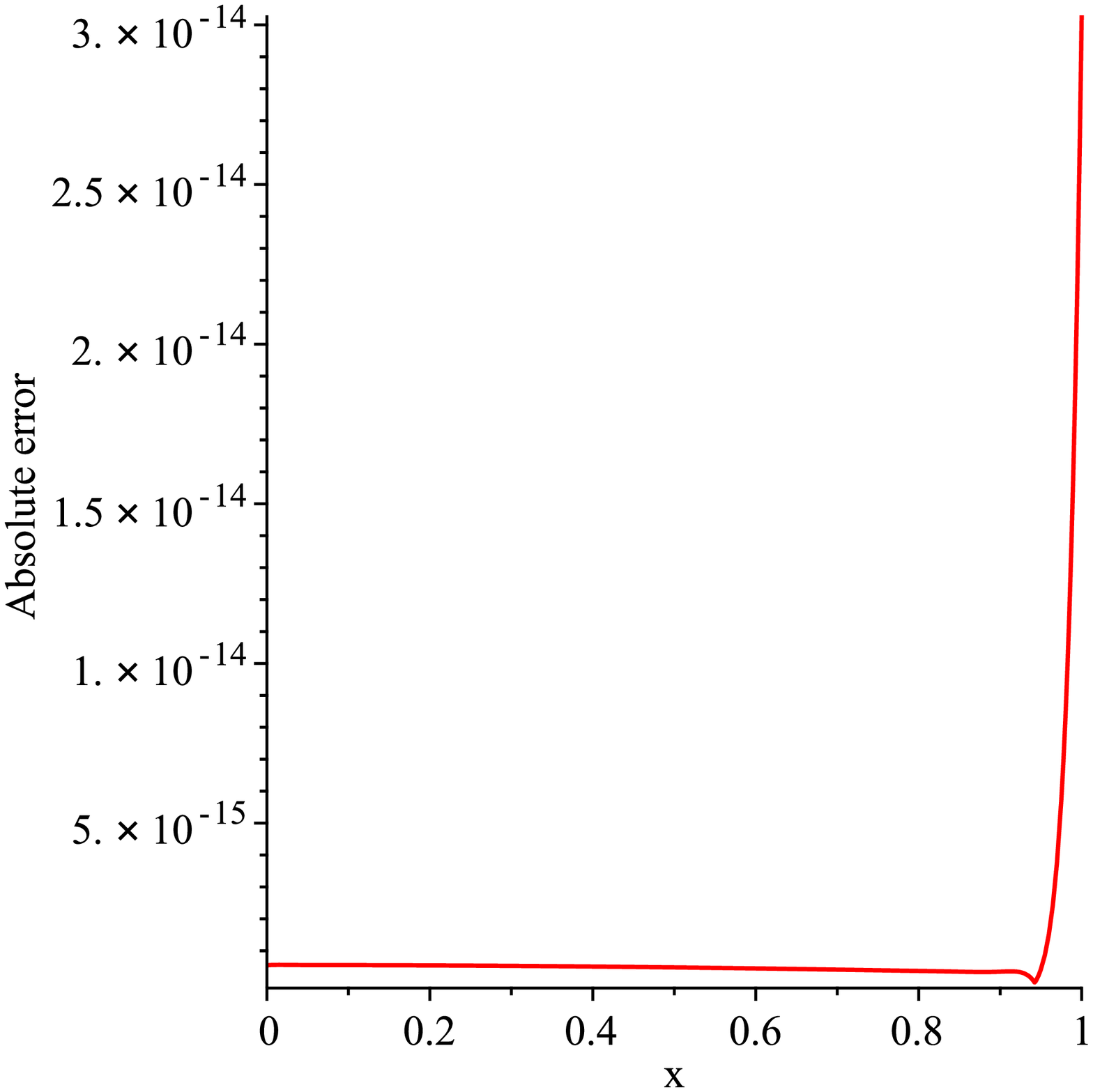}
		\includegraphics[width=2.400in,height=1.55in]{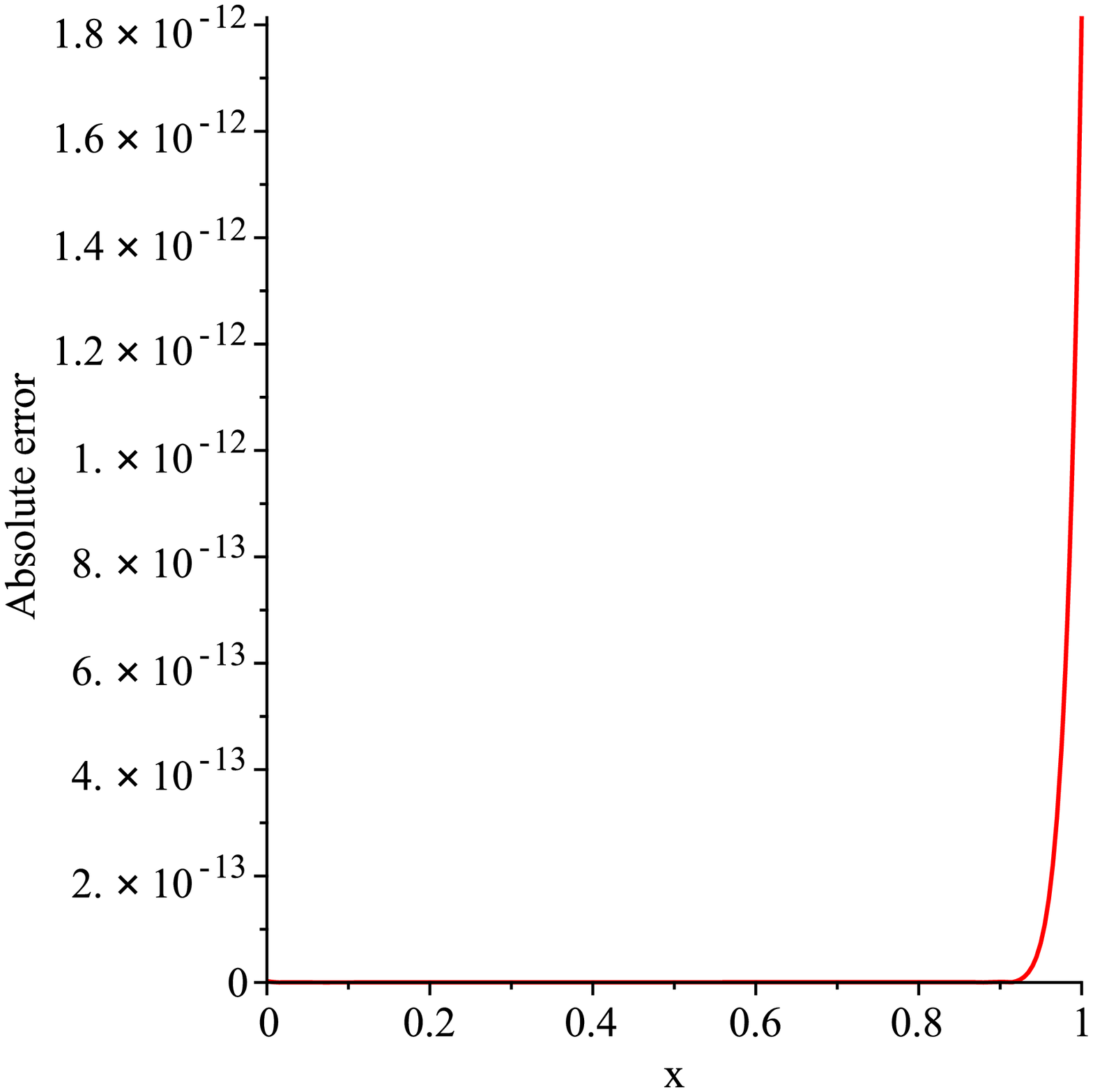}}
	\caption{The plot shows
		exact solution and L-RKM solution of Example 5.1 for $\alpha$=2, $\lambda=1$, $m=20$ and $n=30$ (on
		left), The plot shows absolute error of L-RKM solution of Example 5.1 for
		$\alpha$=2, $\lambda=1$, $m=20$ and $n=30$ (on right).}
	\label{fig1}
\end{figure*}

\begin{figure*}[htbp]
	\centerline{\includegraphics[width=2.400in,height=1.55in]{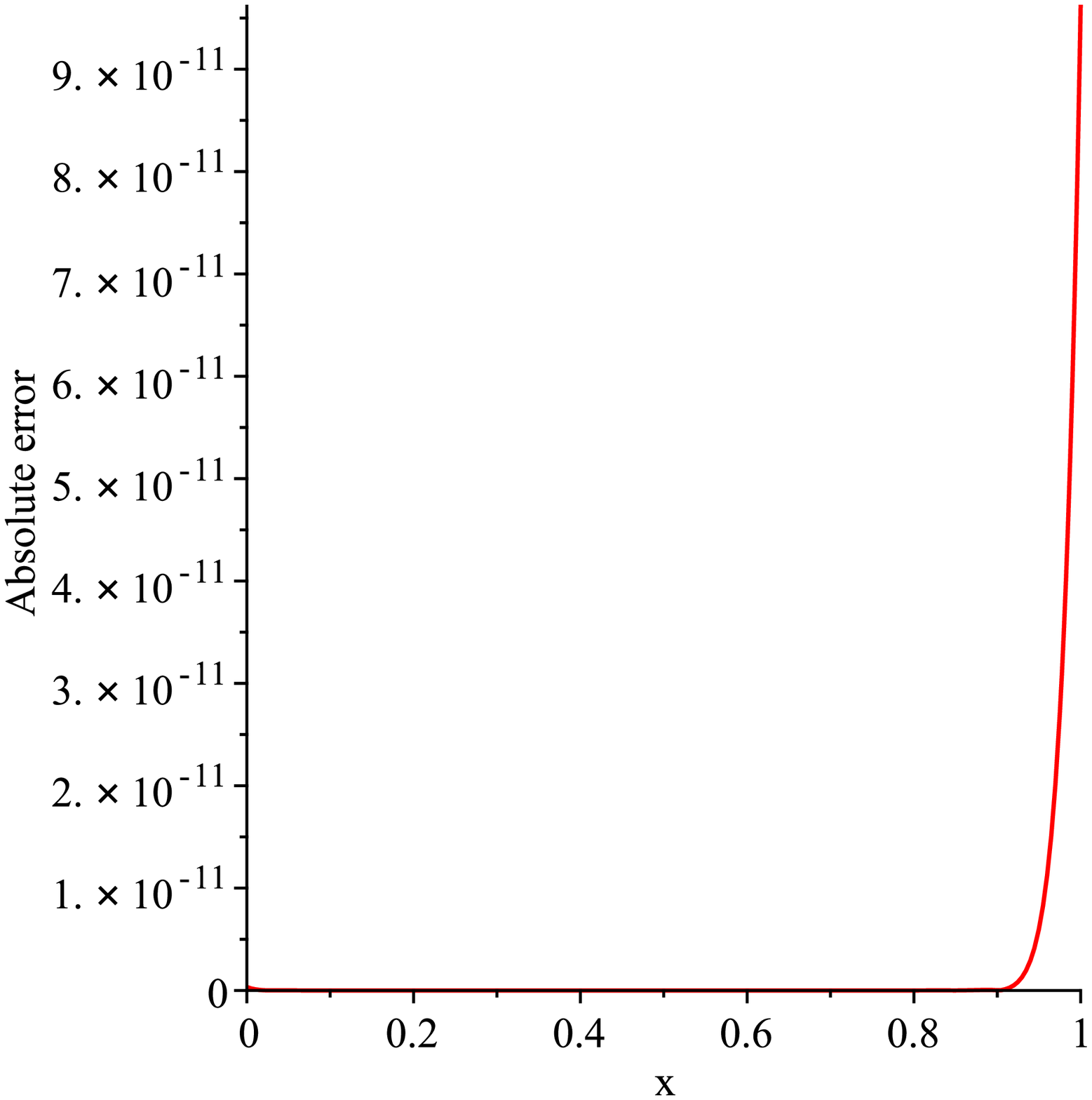}
		\includegraphics[width=2.400in,height=1.55in]{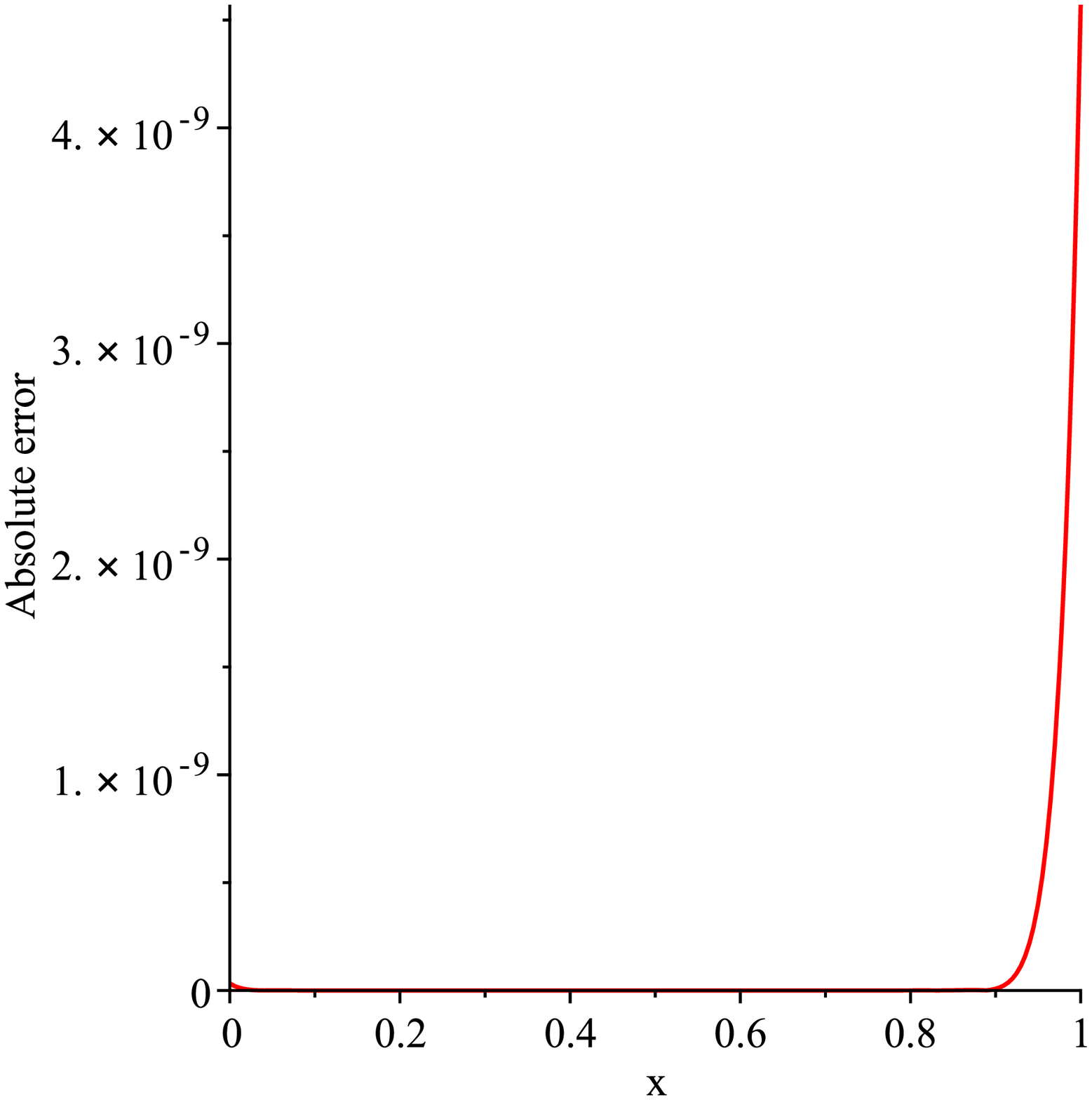}}
	\caption{The plot shows
		absolute error of $y^{\prime}(x)$ for Example 5.1 with $\alpha$=2, $\lambda=1$, $m=20$ and $n=30$ (on
		left), The plot shows absolute error of $y^{\prime\prime}(x)$ for Example 5.1
		with $\alpha$=2, $\lambda=1$, $m=20$ and $n=30$ (on right).}
	\label{fig2}
\end{figure*}

\begin{figure*}[htbp]
	\centerline{\includegraphics[width=2.400in,height=1.55in]{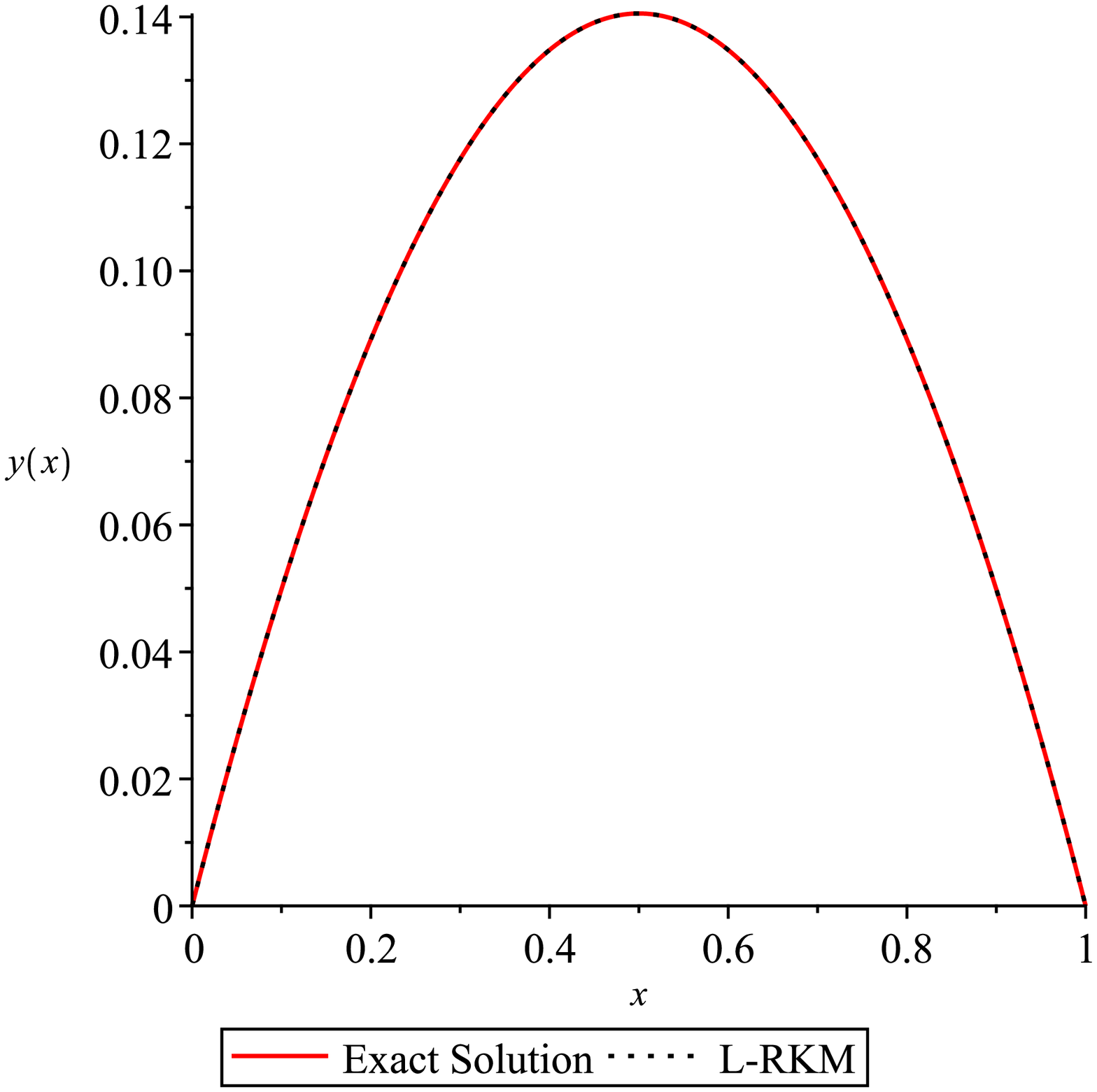}
		\includegraphics[width=2.400in,height=1.55in]{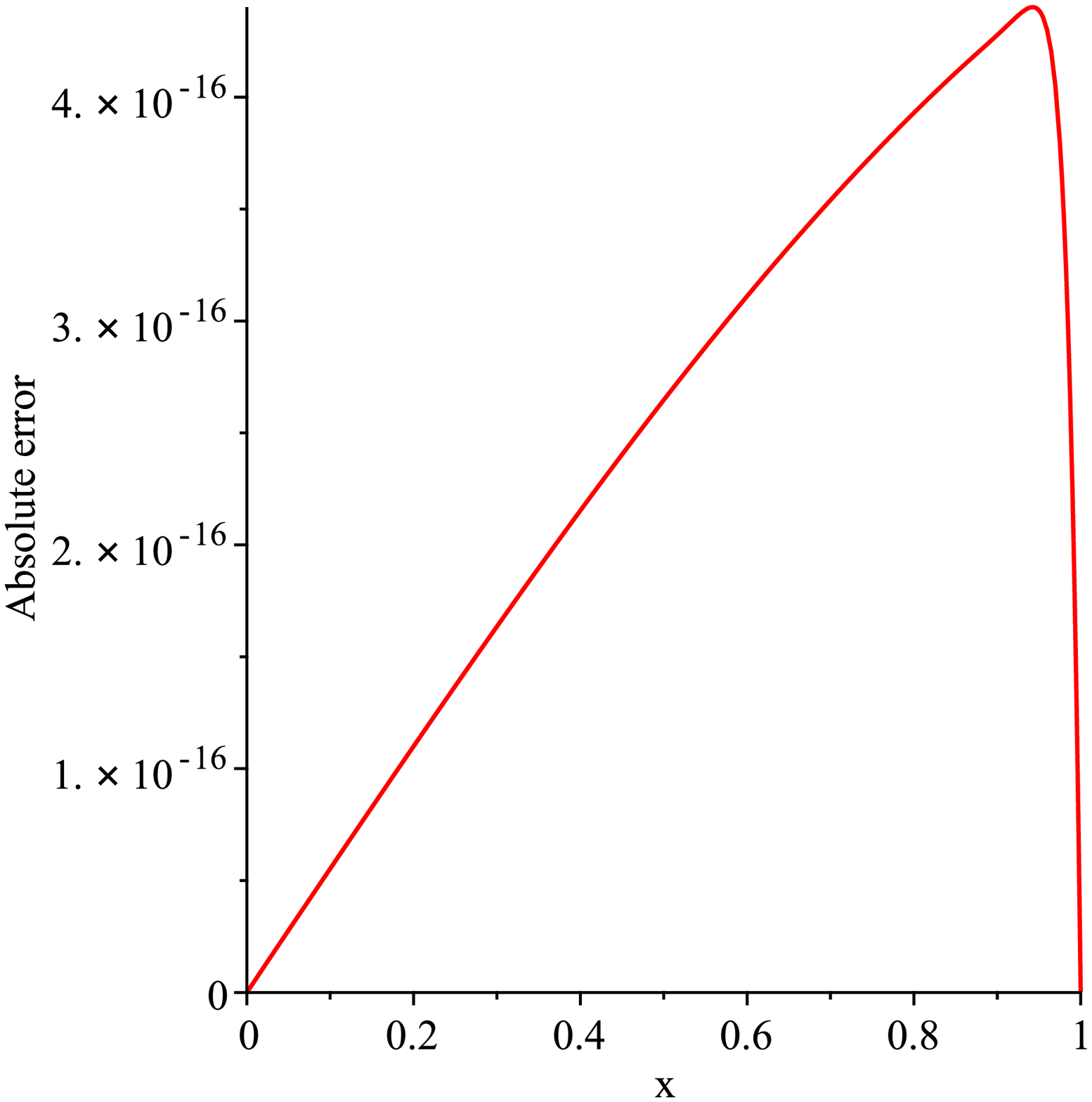}}
	\caption{The plot shows
		absolute error of $y^{\prime\prime\prime}(x)$ for Example 5.1 with $\alpha$=2, $\lambda=1$, $m=20$ and $n=30$ (on
		left), The plot shows absolute error of $y^{(4)}(x)$ for Example 5.1
		with $\alpha$=2, $\lambda=1$, $m=20$ and $n=30$ (on right).}
	\label{fig3}
\end{figure*}

\begin{figure*}[htbp]
	\centerline{\includegraphics[width=2.400in,height=1.55in]{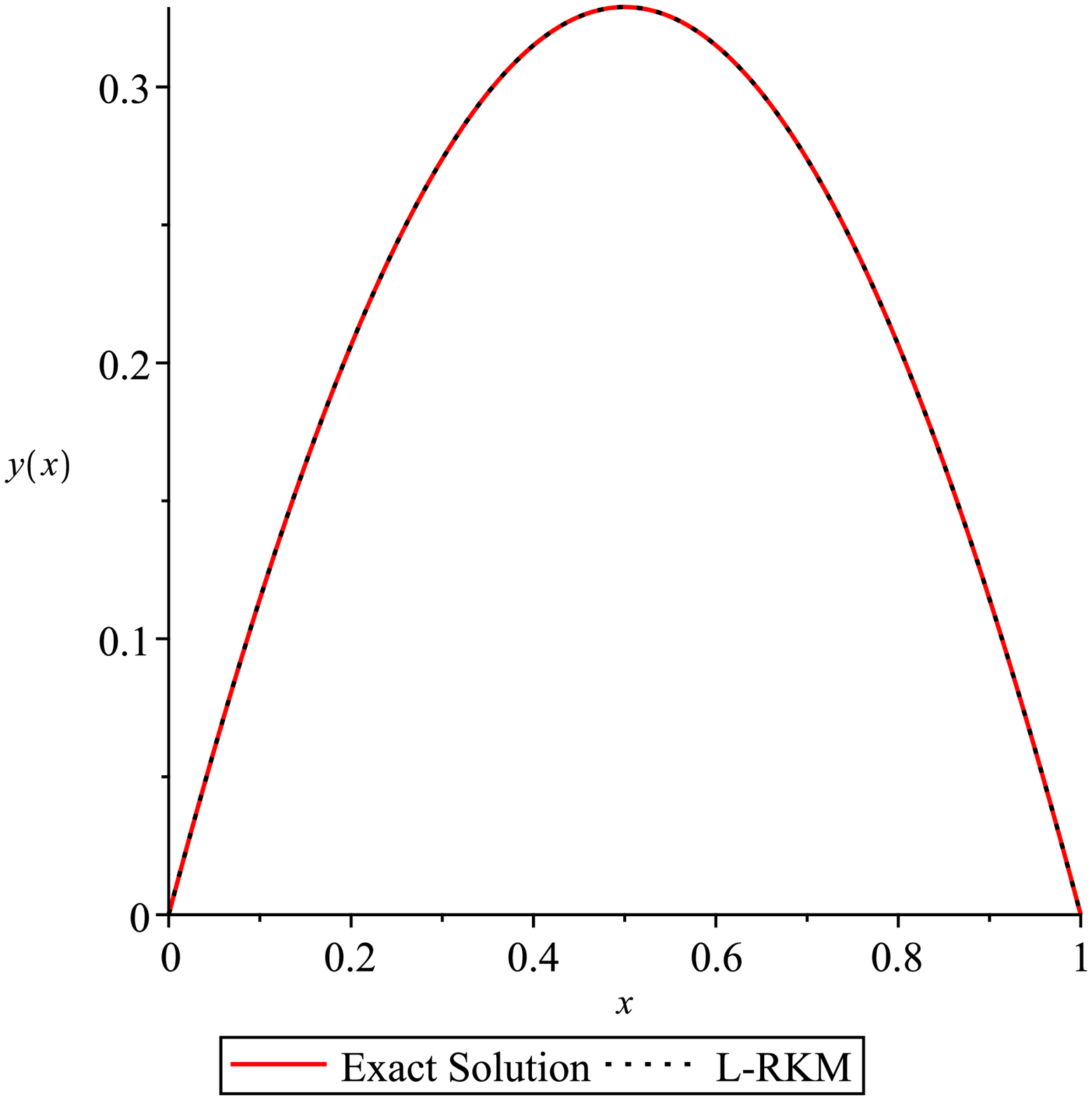}
		\includegraphics[width=2.400in,height=1.55in]{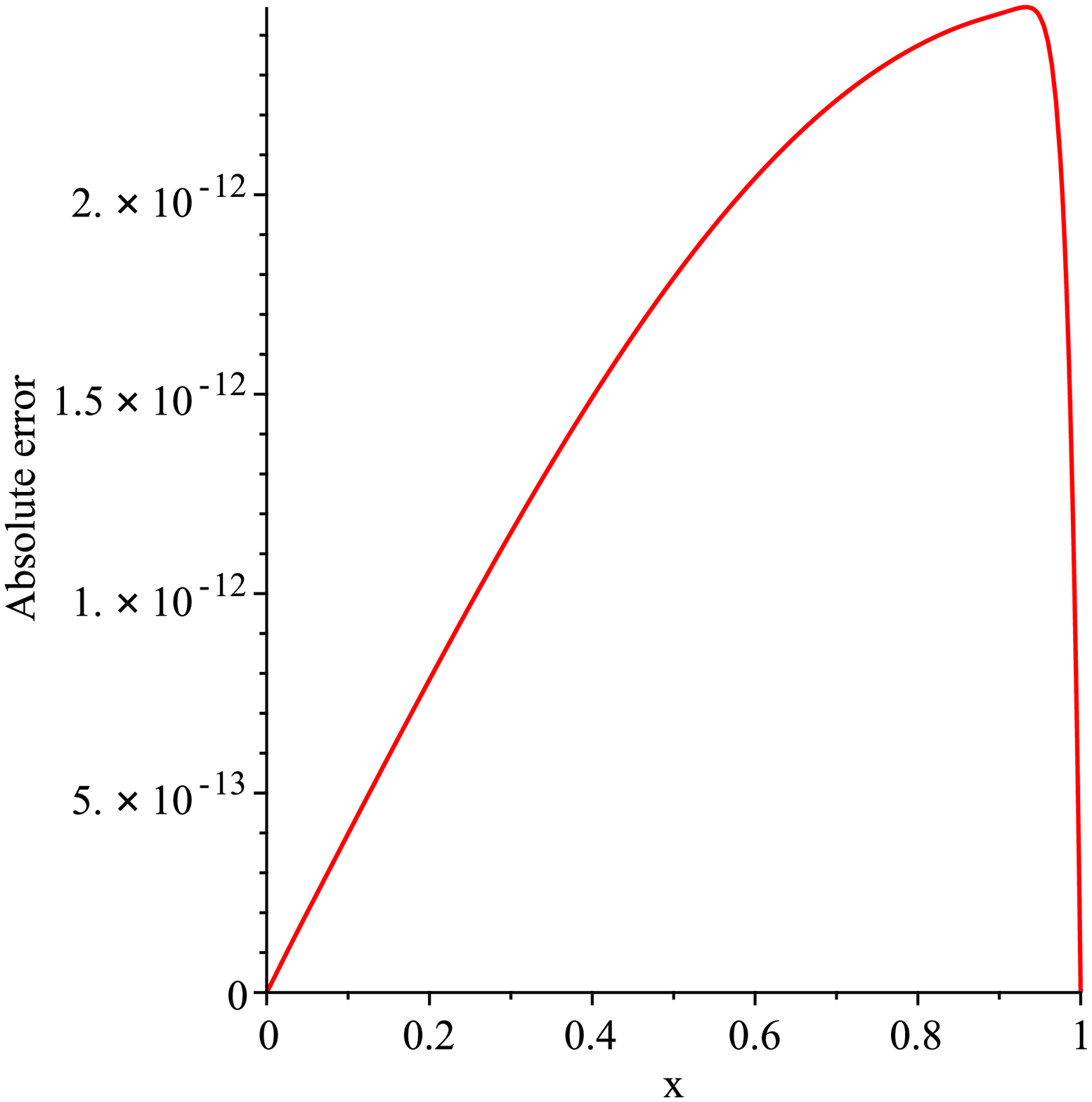}}
	\caption{The plot shows
		exact solution and L-RKM solution of Example 5.1 for $\alpha$=2, $\lambda=2$, $m=20$ and $n=30$ (on
		left), The plot shows absolute error of L-RKM solution of Example 5.1 for
		$\alpha$=2, $\lambda=2$, $m=20$ and $n=30$ (on right).}
	\label{fig4}
\end{figure*}

\begin{figure*}[htbp]
	\centerline{\includegraphics[width=2.400in,height=1.55in]{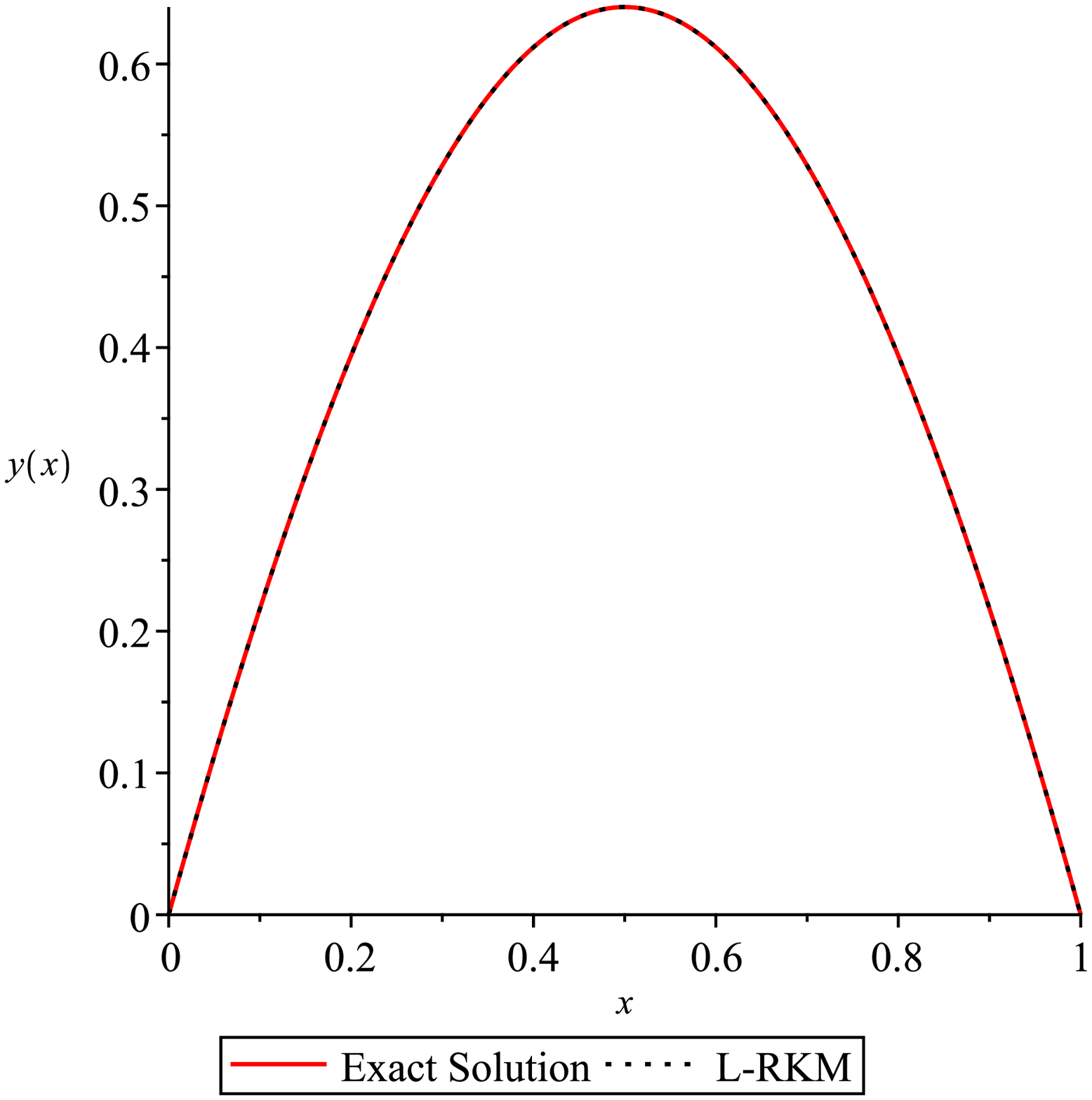}
		\includegraphics[width=2.400in,height=1.55in]{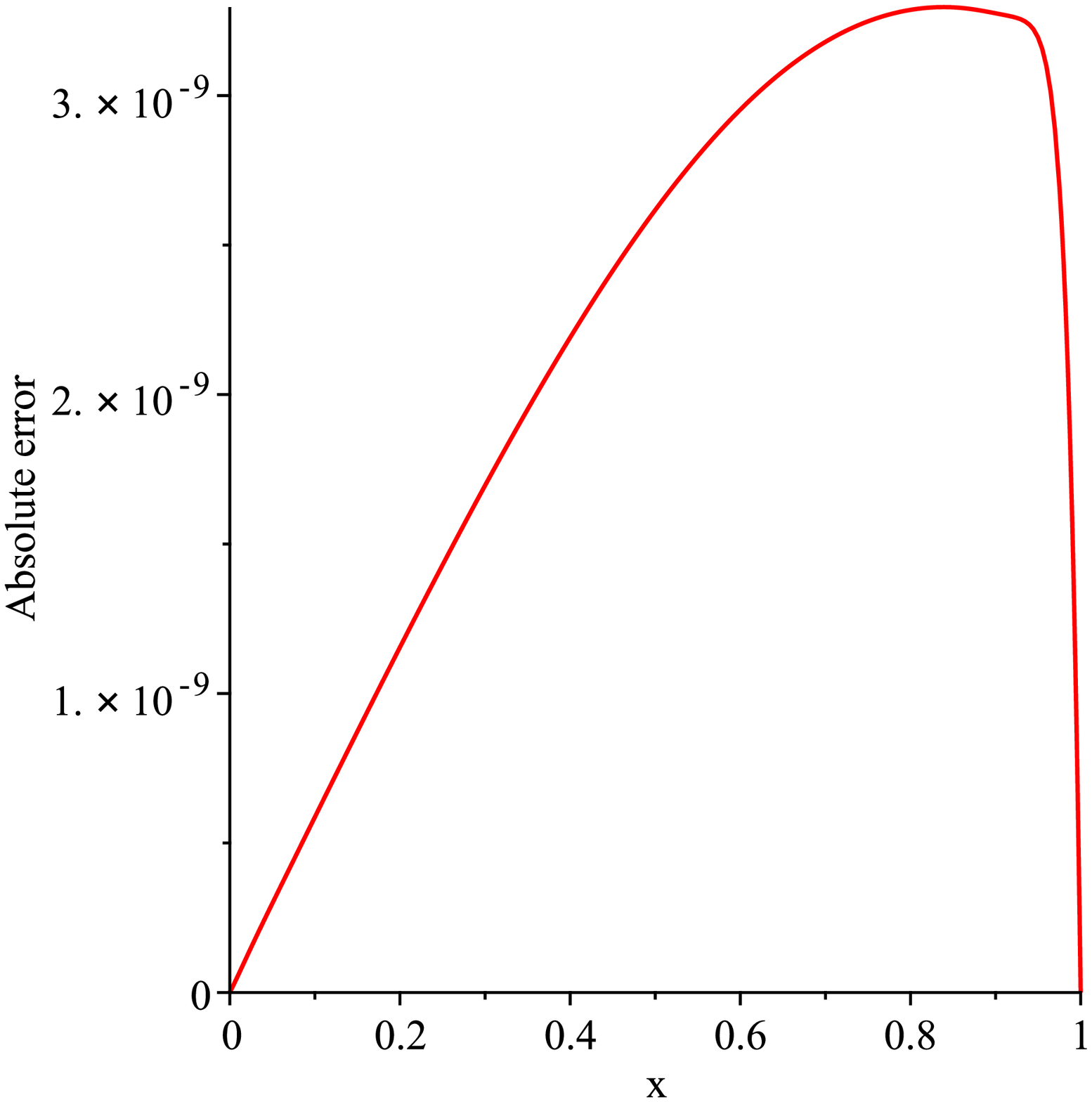}}
	\caption{The plot shows
		exact solution and L-RKM solution of Example 5.1 for $\alpha$=2, $\lambda=3$, $m=20$ and $n=30$ (on left), The plot shows absolute error of L-RKM solution of Example 5.1 for
		$\alpha$=2, $\lambda=3$, $m=20$ and $n=30$ (on right).}
	\label{fig5}
\end{figure*}
\end{document}